\documentclass{amsart}

\usepackage{amsthm,amsmath,amssymb,amscd,amsfonts}

\theoremstyle{plain}
\newtheorem{theorem}{Theorem}[section]
\newtheorem{lemma}[theorem]{Lemma}
\newtheorem{proposition}[theorem]{Proposition}

\theoremstyle{definition}

\newtheorem*{acknowledgement}{Acknowledgement}

\theoremstyle{remark}

\newcommand{\C}{{\mathbb{C}}}

\renewcommand{\P}{{\mathbb{P}}}

\DeclareMathOperator{\genus}{{genus}}
\DeclareMathOperator{\Dom}{{Dom}}
\DeclareMathOperator{\rank}{{rank}}

\DeclareMathOperator{\Chow}{{Chow}} 
\DeclareMathOperator{\id}{{id}}
\DeclareMathOperator{\pic}{Pic^0} 
\begin{document}

\title[Uniformly effective Shafarevich Conjecture]{Uniformly effective Shafarevich Conjecture on families of
hyperbolic curves over a curve with prescribed degeneracy locus}

\author{Gordon Heier}

\address{Ruhr-Universit\"at Bochum\\
Fakult\"at f\"ur Mathematik\\
D-44780 Bochum\\
Germany}

\email{heier@cplx.ruhr-uni-bochum.de}


\begin{abstract}
The paper's main result is an effective uniform bound for the finiteness statement of the Shafarevich Conjecture over function fields. Several results on the projective geometry of curves are established in the course of the proof. These results should be of independent interest. As a corollary, a uniform effective bound for the Mordell Conjecture over function fields is derived via Parshin's trick.
\end{abstract}

\maketitle

\section{A uniform effective solution to the Shafarevich Conjecture}\label{shafsection}
\subsection{Statement of the uniform effective bound}
Let $B$ be a smooth complex projective curve of genus $q\geq 0$. Let
$S\subset B$ be a finite subset of cardinality $s$. The following statement
was conjectured by Shafarevich and proved by Parshin (\cite{P}, case
$S=\emptyset$) and Arakelov (\cite{A}).
\begin{theorem}[\cite{P}, \cite{A}]\label{shaf}
Let $g\geq 2$. Then there are only a finite number of isomorphism classes of
nonisotrivial minimal families $f:X\to B$ of curves of genus $g$ with $X$
smooth such that $f:X\backslash f^{-1}(S)\to B\backslash S$ is a smooth
family.
\end{theorem}
Recall that a family of curves is called {\it isotrivial} if its smooth
fibers are all isomorphic to each other.\par
In the article \cite{C}, Caporaso makes the point that the number of nonisotrivial families in Theorem \ref{shaf} can be bounded by
a uniform constant depending only on $(g,q,s)$. The proof given consists, like ours, of a detailed analysis of the boundedness aspect of the problem. However, the arguments used in \cite{C} are ineffective in nature and differ in essential ways from ours. In fact, we will be able to do without much of the advanced moduli theory used in \cite{C}, since we replace the use of $\overline{M}_g$ with more straightforward algebraic geometric arguments involving Chow varieties, which enable us to argue effectively.\par
Our main theorem is Theorem \ref{effshafconj}. It is fair to say that its proof was inspired by the ``boundedness and rigidity"-type proofs of (effective) finiteness theorems for maps between hyperbolic complex manifolds (see e.g. \cite{MDLM}, \cite{HS}, \cite{Maehara}, \cite{BD}, \cite{TsaiIMRN}, \cite{TsaiCrelle}, \cite{TsaiJAG}, \cite{Guerra}, \cite{heierfinitemaps}). However, our proof is certainly not a straightforward generalization, and we will need to prove several new algebraic geometric results that should be of independent interest.\par
Note that the values of the constants $A,D,Q$ appearing in Theorem \ref{effshafconj} will be defined only after the Theorem has been stated. The expressions that occur are given in the very form that they arise in naturally in the course of our argument. Assuming that any further manipulations are more likely to confound the reader than to help him, we will make no effort to expand the terms below in any way. Obviously, the bound achieved involves iterated exponentials and is most likely far from being sharp. 
\begin{theorem}\label{effshafconj}
Let $g,q\geq 2$. Then the number of isomorphism classes of families $f:X\to B$ described in {\rm Theorem \ref{shaf}} is at most:
\begin{equation*}
(6(q-1)+Q\cdot D)\cdot {5(q-1)A\cdot (6(q-1)+Q\cdot D))\choose 5(q-1)A-1}^{5(q-1)A((6(q-1)+Q\cdot D)^2+1)}.
\end{equation*}
In the two remaining cases $q=0,1$, the bound we seek can be taken to be
the above bound with $q$ replaced by $2$ and $s$
replaced by $2s$, multiplied by $\mathcal S(g)$.
\end{theorem}
We start the definition of the constants by defining $\mathcal S(g)$. It can be taken to be any bound on the number of holomorphic maps with a given smooth compact complex curve of genus $g\geq 2$ as domain and any smooth compact complex curve of genus at least $2$ as target. The fact that the number in question is finite is the classical Theorem of de Franchis-Severi and an effective bound is relatively easy to obtain with the techniques from the above cited papers. For simplicity, we adopt the result of \cite{HS} and let
\begin{equation*}\label{HSbound}
\mathcal S(g) := 42(g-1)(\frac 1 2 (2\sqrt{6}(g-1)+1)^{2+2g^2}g^2(g-1)(\sqrt{2})^{g(g-1)}+1).
\end{equation*}\par
If we let
\begin{equation*}
l=4\cdot1250(gq+s)-3,
\end{equation*}
then the constants $Q,D,A$ are determined by
\begin{eqnarray*}
Q&=&l\cdot5(2g-2)\left({5\cdot 1250(gq+s)-3\choose 1250(gq+s)}-1\right)\\
&&\cdot {4(4\cdot1250(gq+s)-3)\cdot5(2g-2)+{5\cdot 1250(gq+s)-3\choose
1250(gq+s)}-1\choose {5\cdot 1250(gq+s)-3\choose 1250(gq+s)}-1},
\end{eqnarray*}
\begin{equation*}
D=\left({l\cdot 5(2g-2)+2\choose 2}-1\right) \cdot l^2\cdot 5(2g-2)\cdot 500((l\cdot 5(2g-2))^2+1)(gq+s)k
\end{equation*}
and
\begin{equation*}
A={ld+2\choose 2}^Q.
\end{equation*}\par
By thinking of $X$ as a Riemann surface for the moment, one is intuitively
led to expect that the number of families in Theorem \ref{shaf} is actually
zero if $q=0$ and $s\leq 2$ or $q=1$ and $s=0$. This is indeed true, and a
proof of this fact can be found in \cite{ViehwegICTP}. We also refer to
\cite{ViehwegICTP} (and to \cite{Viehwegbook}) for an account of the more
recent history of the theory of families of manifolds and the related
questions about the positivity of direct image sheaves.
\subsection{Sketch of the proof}
The following paragraphs contain an outline of the strategy used to prove Theorem \ref{effshafconj}. First, we assume $q\geq 2$.\par
Let $f:X\to B$ be one of the families under discussion. Since $q\geq 2$, $X$
is a manifold of general type and there exist a finite number of rational
$(-2)$-curves $C_i$ $(i\in I)$ with $\bigcup_{i\in I}C_i\subset f^{-1}(S)$
such that, by a theorem of Bombieri, the $5$-canonical map
\begin{equation*}
\varphi_{|5K_X|}:X\to \P(H^0(X,5K_X))=:\P^{m_X}
\end{equation*}
is an embedding on the complement of $\bigcup_{i\in I}C_i$ and contracts the
$C_i$ to rational double point singularities. In \cite[Proposition 2]{P} it
is stated that $m_X\leq 50\cdot 5^2(gq+s)=1250(gq+s)=:m$ as a consequence of
the Riemann-Roch Theorem. We can assume (after linear inclusions) that $m_X=m$ for all families.\par
Since the degree of the divisor $5K_X$ on the smooth fibers $X_b$ is equal to $d:=5(2g-2)$
(independent of $b$), there corresponds to every family $f:X\to B$ a
canonical morphism $\psi_X^0:B\backslash S\to\Chow_{1,d}(\P^m)$ given by
$b\mapsto [\varphi_{|5K_X|}(X_b)]$. We follow the standard convention that $\Chow_{1,d}(\P^m)$ denotes the Chow variety of $1$-dimensional cycles of degree $d$ in $\P^m$. An introduction to the theory of Chow varieties can be found for example in \cite{TsaiJAG}, or in many textbooks on algebraic geometry. Since $B$ is smooth and
$\Chow_{1,d}(\P^m)$ is projective, there exists a unique extension
$\psi_X:B\to\Chow_{1,d}(\P^m)$ that coincides with $\psi_X^0$ on
$B\backslash S$.\par
Clearly, no two nonisomorphic families with the properties described in
Theorem \ref{shaf} can correspond to the same isomorphism class of morphisms $\psi:B\to\Chow_{1,d}(\P^m)$.
Thus, we are left with bounding the cardinality of the set of isomorphism classes of such maps
$\psi_X$, which is achieved by bounding the
degree of the graphs $\Gamma_{\psi_X}\subset B\times\Chow_{1,d}(\P^m)$ and
then applying a rigidity argument to $\Gamma_{\psi_X}$ in conjunction with an
estimate on the complexity of a certain Chow variety of cycles in $B\times\Chow_{1,d}(\P^m)$. This finishes the proof in the
case $q\geq 2$.\par
The remaining cases $q=0,1$ are simply handled by executing a degree $2$
base change to a curve of genus $2$ and thus returning to the situation
dealt with previously.
\subsection{A degree bound for the image of the moduli map}
The main difficulty in bounding $\deg \Gamma_{\psi_X}$ lies in bounding the
degree of $\psi_X(B)\subset\Chow_{1,d}(\P^m)$. The key to achieving this will lie in the construction of a very natural yet somewhat non-standard embedding of the Chow variety in question into a product of projective spaces. In this product, degree is measured with respect to the Segre embedding, making it sufficient to be able to control the degree of $\psi_X(B)\subset\Chow_{1,d}(\P^m)$ under the component maps, which will be written as
\begin{equation*}
\Phi_\nu:\Chow_{1,d}(\P^m)\to\Chow_{1,ld}(\P^2_\nu)\cong\P^{{ld+2\choose
2}-1}_\nu.
\end{equation*}\par
The construction of the aforementioned embedding is one of the main results in this paper.
\subsubsection{A generalization of a result of Clemens}
The following Proposition \ref{Clgen} constitutes a generalization of Clemens' result in
\cite[Theorem 1.1]{Cl}. However, our approach to the proof is quite
different from Clemens' original approach; in fact, our proof is more in the
spirit of Ein's papers \cite{EinI} and \cite{EinII} and Voisin's paper
\cite{Voisin} (to which a substantial correction had to be
published in \cite{Voisincorrection}). Notwithstanding, our proof contains a
key new element. Namely, we will be using the technique of explicit
constructions of meromorphic vector fields of low pole order which was first
used in the present manner by Siu in \cite{SiuICM}. This idea makes the proof less technical and should allow for even further generalizations. The author will pursue this avenue in future research. \par
The following Proposition will be used in Section \ref{twistedproj} to show the
well-definedness of the map $\Phi$ which we are going to construct.
\begin{proposition}\label{Clgen}
Let $f_0,f_1,f_2$ be a generic triple of homogeneous polynomials of degree
$l\geq \frac{2g-2}{\delta}+(4m-4)$ in the $m+1$ homogeneous coordinates of
$\P^m$. Then, for any irreducible complex curve $C$ of genus $g$ and degree
$\delta$ in $\P^m$, the topological closure in $\P^2$ of the image of
$C\backslash\{[Z]\in \P^m:f_0(Z)=f_1(Z)=f_2(Z)=0\}$ under the holomorphic
map
\begin{equation*}
\phi_{f_0,f_1,f_2}:\P^m\backslash\{[Z]\in \P^m:f_0(Z)=f_1(Z)=f_2(Z)=0\}\to
\P^2
\end{equation*}
defined by $[Z]\mapsto [f_0(Z),f_1(Z),f_2(Z)]$ is a curve (necessarily of
degree $l\delta$, when multiplicity is counted) in $\P^2$.
\end{proposition}
\begin{proof}
What we need to show is that the image of $C$ under $\phi_{f_0,f_1,f_2}$ is
$1$-dimensional for a generic choice of $f_0,f_1,f_2$. First, we can assume w.l.o.g. that none of the $f_i$ vanish identically on $C$, since this would represent a trivial special case of the subsequent argument. Now, given that none of
the $f_i$ vanish identically on $C$, the image of $C$ under $\phi_{f_0,f_1,f_2}$ is not $1$-dimensional if and only if
\begin{equation}\label{topoint}\exists
\alpha,\beta\in\C\backslash\{0\}:(f_0-\alpha f_1)|_C\equiv 0 \text{ and }
(f_0-\beta f_2)|_C\equiv 0.
\end{equation}
We shall now argue as follows. Set
\begin{equation*}
f^{(a,b,c)}:=af_0+bf_1+cf_2
\end{equation*}
with $a,b,c\in \C$. Moreover, set $\alpha':=b+\alpha(a+\frac c \beta)$. Now note the equality:
\begin{eqnarray}
&&f^{(a,b,c)}-\alpha'f_1\nonumber\\
&=&a(f_0-\alpha f_1)+c(f_2-\frac \alpha \beta f_1)\label{thepone}.
\end{eqnarray}
If property \eqref{topoint} holds, the expression in \eqref{thepone} clearly vanishes when restricted to $C$. Furthermore, it is evident from \eqref{thepone} that $f^{(a,b,c)}-\alpha'f_1$ represents a linear $\P^1$ represented by homogeneous coordinates $[a,c]$ and contained in the projective $2$-plane spanned by $[f_0],[f_1],[f_2]$.\par
In other words, if the closure of the image of $C$ under
$\phi_{f_0,f_1,f_2}$ is not $1$-dimensional, then the projective 2-plane
spanned by $[f_0],[f_1],[f_2]$ in $\P^{{l+m\choose m}-1}$ contains a linear
$\P^1$ representing equivalence classes of homogeneous polynomials that
vanish identically on $C$. However, we will now prove that $\P^{{l+m\choose
m}-1}$ does not contain a linear $S=\P^{{l+m\choose m}-2}$ of homogeneous
polynomials that vanish identically on $C$. In fact, we will show that this is not possible even if we allow different curves $C$ for different polynomials. Once we have proved this, we can choose $[f_0],[f_1],[f_2]$ to be
so general in $\P^{{l+m\choose m}-1}$ that the $2$-plane spanned by them
does not contain a linear $\P^1$ representing equivalence classes of
homogeneous polynomials that vanish identically on $C$. The proof of the
Proposition will then be completed.\par
Let $S$ be a linear hyperplane of $\P^{{l+m\choose m}-1}$. In homogeneous coordinate
systems $[Z_i]$ ($i=0,\ldots,m$) on $\P^m$ and $[A_\mu]$ ($|\mu|=l$) on
$\P^{{l+m\choose m}-1}$, let $S$ be given by
\begin{equation*}\{\sum_{|\mu=(\mu_0,\ldots,\mu_m)|=l} b_\mu A_\mu=0\}
\end{equation*}
for $b_\mu\in\C$.
Moreover, let
\begin{equation*}\mathfrak X :=\{([Z_0,\ldots,Z_m],[A_\mu])\in\P^m\times
S|\sum_{|\mu|=l} A_\mu Z^\mu=0\}.
\end{equation*}
Note that $\mathfrak X$ is a smooth hypersurface in $\P^m\times S$. In
particular,
\begin{equation*}
\dim_\C(\mathfrak X)=m+{l+m\choose m}-3.
\end{equation*}\par
We claim that, at every point $P$ of $\mathfrak X$, the vector bundle
$T\mathfrak X\otimes p_1^*\mathcal O_{\P^m}(2)\otimes p_2^*\mathcal
O_{S}(1)$ is generated by global sections. ($p_1:\mathfrak X\to \P^m$ and $p_2:\mathfrak X\to S$ denote the two canonical projections from $\mathfrak X$). \par
To see this, first note that it suffices to prove the claim at a general
$P$, since there exists a finite number of homogeneous coordinate systems on $\P^m$ and
$\P^{{l+m\choose m}-1}$ such that for every point of $\mathfrak X$,
there is one pair of coordinates for which that point is general in the
sense given below.\par
We continue the proof of our claim by writing down the following explicit
meromorphic vector fields on $\mathfrak X$.
\begin{equation}\label{vf1}
Z_0\frac{\partial}{\partial
Z_0}-\sum_{|\mu|=l}\mu_0A_\mu\frac{\partial}{\partial A_\mu
}-\frac{\sum_{|\mu|=l}\mu_0A_\mu b_\mu }{\sum_{|\mu|=l}\mu_i A_\mu b_\mu
}(Z_i\frac{\partial}{\partial
Z_i}-\sum_{|\mu|=l}\mu_iA_\mu\frac{\partial}{\partial A_\mu })
\end{equation}
for $i=1,\ldots,m-1$;
\begin{eqnarray}\label{vf2}
&&\frac{Z_i}{Z_0}A_{\lambda +e_j}\frac{\partial}{\partial A_{\lambda
+e_j}}-\frac{Z_j}{Z_0}A_{\lambda +e_i}\frac{\partial}{\partial A_{\lambda
+e_i}}\\
&+&\left(\frac{Z_jA_{\lambda +e_i}b_{\lambda +e_i}-Z_iA_{\lambda
+e_j}b_{\lambda +e_j}}{Z_2 A_{\lambda^{(0)} +e_1}b_{\lambda^{(0)} +e_1}-Z_1
A_{\lambda^{(0)} +e_2}b_{\lambda^{(0)} +e_2}}\right.\nonumber\\
&&\cdot\left.(\frac{Z_2}{Z_0}A_{\lambda^{(0)} +e_1}\frac{\partial}{\partial
A_{\lambda^{(0)} +e_1}}-\frac{Z_1}{Z_0}A_{\lambda^{(0)}
+e_2}\frac{\partial}{\partial A_{\lambda^{(0)} +e_2}})\right),\nonumber
\end{eqnarray}
where $|\lambda|=l-1;0\leq i< j\leq m$ and $\lambda^{(0)}$ is any fixed
index.\par
The first thing to notice is that these vector fields are all tangent to
$\mathfrak X$ (apply them to the defining equations $\sum_{|\mu|=l} b_\mu
A_\mu$ and $\sum_{|\mu|=l} A_\mu Z^\mu$ and you will get zero). Second,
their pole order is no more than $2$ in the $Z$-direction and no more than
$1$ in the $A$-direction.\par
Third, we can assume w.l.o.g.\ that in our coordinate system, the
coefficients of the above vector fields will all be nonzero finite numbers
at $P$. It is obvious that the vector fields in \eqref{vf1} generate a
linear space of dimension $m-1$ in $T_P(\mathfrak X)$. After a lexicographical ordering of the indices $\lambda$, a simple count yields that
the vector fields in \eqref{vf2} generate a linear space of dimension
${l+m\choose m}-2$ in $T_P(\mathfrak X)$ that shares only the zero vector with the span
of the vectors of \eqref{vf1}. Therefore, the dimension of the span of the
union of the vector fields in \eqref{vf1} and \eqref{vf2} in $T_P(\mathfrak
X)$ is
\begin{equation*}
m-1+{l+m\choose m}-2=m+{l+m\choose m}-3,
\end{equation*}
which is what we needed to show to prove our claim about the generation of $T\mathfrak X\otimes p_1^*\mathcal O_{\P^m}(2)\otimes p_2^*\mathcal
O_{S}(1)$.\par
Now assume, in order to derive a contradiction, that indeed every
homogeneous polynomial represented by a point $s \in S$ vanishes on some
curve $C_s$ of degree $\delta$. Then there exists a commutative diagram of
families of curves
\begin{equation*}
\begin{CD}
\tilde {\mathfrak C} @>>>\mathfrak C@>>>\mathfrak X\\
@VVV @VVV @VVp_2V\\
\tilde S @>\eta>> S@>=>>S\\
\end{CD}
\end{equation*}
such that $\tilde {\mathfrak C}\to \tilde S$ is a family with smooth total
space and $\eta:\tilde S\to S$ is a covering map over a Zariski open
subset of $S$ such that $\tilde C_{\tilde s}$ is a normalization of
$C_{\eta(\tilde s)}$ for generic $\tilde s\in \tilde S$.
Note that, in the remainder of our proof, $S$ can be replaced with any small
open subset of itself, so we can abuse notation and assume w.l.o.g.\ that
$\tilde {\mathfrak C}\to \tilde S$ and $ \mathfrak C\to S$ are identical.\par
We have the short exact sequence (with $L:=p_1^*\mathcal O_{\P^m}(2)\otimes
p_2^*\mathcal O_{S}(1)$)
\begin{eqnarray*}
0\to T\mathfrak C\otimes L|_\mathfrak C \to T\mathfrak X\otimes L|_\mathfrak
C\to N_{\mathfrak C,\mathfrak X}\otimes L|_\mathfrak C\to 0.
\end{eqnarray*}
Since $T\mathfrak X\otimes L|_\mathfrak C$ is generated by global sections,
the first Chern class of the restriction of $T\mathfrak X\otimes L|_\mathfrak
C$ to a generic curve $C_s$ is nonnegative and so is the first Chern class
of the restriction of the quotient bundle $N_{\mathfrak C,\mathfrak
X}\otimes L|_\mathfrak C$ to $C_s$. For a generic $s\in S$, we have
\begin{equation*}
N_{\mathfrak C,\mathfrak X}\otimes L|_{C_s}\cong N_{C_s,X_s}\otimes L|_{C_s},
\end{equation*}
giving
\begin{equation}\label{chernpos}
c_1(N_{C_s,X_s}\otimes L|_{C_s})\geq 0.
\end{equation}
For $C_s$, there is the short exact sequence
\begin{equation*}
0\to TC_s\otimes L|_{C_s} \to TX_s\otimes L|_{C_s}\to N_{C_s,X_s}\otimes L|_{C_s}\to
0.
\end{equation*}
Because of \eqref{chernpos}, we get
\begin{equation*}
c_1( TX_s\otimes L|_{C_s})\geq c_1(TC_s\otimes L|_{C_s}),
\end{equation*}
or equivalently
\begin{eqnarray*}
&&(m-1)3\delta+c_1(TX_s|C_s)\geq 2-2g+3\delta\\
&\Leftrightarrow&(m-1)3\delta+(m+1)\delta-l\delta\geq 2-2g+3\delta\\
&\Leftrightarrow&(m-1)3+(m+1)-l\geq \frac{2-2g}{\delta}+3\\
&\Leftrightarrow&4m-5-l\geq \frac{2-2g}{\delta}\\
&\Leftrightarrow&l\leq\frac{2g-2}{\delta}+(4m-5),\\
\end{eqnarray*}
which contradicts the assumption on $l$.
\end{proof}
\subsubsection{Construction of the key embedding}\label{twistedproj}
The following lemma is the key to the injectivity statement which we will
prove at the end of this section.
\begin{lemma}\label{injlin}
Let $\delta_0\geq 1$ and $M\geq 3$ be integers. Then there exist linear
projection maps
\begin{equation*}
\pi_\nu:\P^M\backslash\P^{M-3}_\nu\to\P^2_\nu \quad
(\nu=1,\ldots,Q:=\delta_0M{4\delta_0+M\choose M})
\end{equation*}
such that the rational map
\begin{equation*}\Chow_{1,\delta_0}(\P^M)\rightharpoonup\prod_{\nu=1}^{Q}\Chow_{1,\delta_0}(\P^2_\nu)
\end{equation*}
given by $[C]\mapsto ([\overline{\pi_\nu(C\cap \Dom(\pi_\nu))}])_{\nu=1,\ldots,Q}$ (where the image curves are taken to have appropriate multiplicities) is
injective on its maximal set of definition, which is the set of those
$[C]\in \Chow_{1,\delta_0}(\P^M)$ such that no irreducible component of $C$
is contained in any of the light sources $\P^{M-3}_\nu$ or is mapped to a point by any of the $\pi_\nu$.
\end{lemma}
\begin{proof}
We first describe the construction of the $Q$ projections required. Let
$H_r$ ($r=1,\ldots,{4\delta_0+M\choose M}$) be ${4\delta_0+M\choose M}$
points in general position in the dual space ${\P^M}^*$. Formulae for the degree of dual varieties are well-known, and we simply quote the one given in \cite[Lemma 1.2]{HS} to establish that the degree of the dual variety of a curve of
degree $2\delta_0$ is no more than $4\delta_0$. Thus, there is at least one
index $\rho$ such that $H_\rho$ is not contained in the dual of $C\cup
\tilde C$ if $C,\tilde C$ are two curves of degree $\delta_0$. In other
words, the linear hyperplane $H_\rho$ intersects $C$ and $\tilde C$
everywhere transversally.\par
Next, for every $r$, choose homogeneous coordinates $[X_0,\ldots,X_M]$ on
$\P^M$ such that $H_r$ is given by $\{X_0=0\}$. For every $r$ and for every $\beta\in \{1,\ldots,\delta_0 M\}$, let $(a_{(r,\beta),\alpha})_{\alpha=1,\ldots,M}$ be an $M$-tuple of complex numbers such that the following nondegeneracy condition is satisfied. For pairwise distinct $\beta^{(1)},\ldots,\beta^{(M)}$:
\begin{equation}\label{nondeg}
\det(a_{(r,\beta^{(\eta)}),\alpha})^{\alpha=1,\ldots,M}_{\eta=1,\ldots,M}\neq 0.
\end{equation}\par
We now define the projection $\pi_{(r,\beta)}$ in the coordinates pertaining
to $H_r$ as
\begin{equation*}
[X_0,\ldots,X_M]\mapsto [X_0,\sum_{\alpha=1}^{M}a_{(r,\beta),\alpha}\cdot X_\alpha,X_2]\ \ \text{(wherever defined)}.
\end{equation*}
Finally, we take $\nu=1,\ldots,\delta_0M{4\delta_0+M\choose M}$ to be a
parameter that counts the tuples $(r,\beta)$ in an arbitrary way.\par
Having defined the projections, let $C$ and $\tilde C$ be two curves of
degree $\delta_0$ in $\P^M$ such that no irreducible component of them is
contained in any of the $\P^{M-3}_\nu$ and such that no irreducible component is mapped to a point by any of the $\pi_\nu$. What we need to prove is that
$C=\tilde C$ if and only if $\pi_\nu(C)=\pi_\nu(\tilde C)$ for all
$\nu$.\par
The ``only if'' part being trivial, we prove the ``if'' part. To this end,
let $\rho$ be such that $H_\rho$ intersects $C\cup\tilde C$ everywhere
transversally. Let $U_j$ ($j\in J$, $\#J=\delta_0$) (resp.~$\tilde U_{\tilde
j}$ ($\tilde j\in \tilde J$, $\#\tilde J=\delta_0$)) be small neighborhoods
in $C$ (resp.~$\tilde C$) around the respective $\delta_0$ points of
intersection with $H_\rho$. For all $j,\tilde j$, there exist local
holomorphic functions $g_\alpha^{(j)}, \tilde g_\alpha^{(\tilde j)}$
($\alpha=1,\ldots,M$) on a neighborhood of $0\in\C$ such that
\begin{equation*} U_j=[x,g_1^{(j)}(x),\ldots,g_{M}^{(j)}(x)]\ \text{and}\ \tilde U_{\tilde j}=[x,\tilde g_1^{(\tilde j)}(x),\ldots,\tilde
g_{M}^{(\tilde j)}(x)].
\end{equation*}\par
The fact that $\pi_\nu(C)=\pi_\nu(\tilde C)$ for all $\nu$ implies that
$\forall j\in J\ \forall \beta=1,\ldots,\delta_0M\ \exists \tilde
j(j,\beta)$ such that
\begin{eqnarray*}
[x,\sum_{\alpha=1}^{M}a_{(r,\beta),\alpha}\cdot g_\alpha^{(j)}(x),g_2^{(j)}(x)] = [x,\sum_{\alpha=1}^{M}a_{(r,\beta),\alpha}\cdot\tilde
g_\alpha^{(\tilde j)}(x),\tilde g_2^{(\tilde j)}(x)].
\end{eqnarray*}
This implies that
\begin{equation}\label{aaa}
\forall j\in J \ \forall \beta=1,\ldots,\delta_0M\ \exists \tilde
j(j,\beta):
\sum_{\alpha=1}^{M}a_{(r,\beta),\alpha}\cdot (g_\alpha^{(j)}(x)-\tilde
g_\alpha^{(\tilde j)}(x))\equiv 0
\end{equation}
and, by symmetry,
\begin{equation}\label{bbb}
\forall \tilde j\in\tilde J \ \forall \beta=1,\ldots,\delta_0M\ \exists
j(\tilde j,\beta):
\sum_{\alpha=1}^{M}a_{(r,\beta),\alpha}\cdot (g_\alpha^{(j)}(x)-\tilde
g_\alpha^{(\tilde j)}(x))\equiv 0.
\end{equation}
Since $\#\tilde J=\delta_0$ and $\beta\in\{1,\ldots,\delta_0M\}$, it follows from
\eqref{aaa} by a simple pigeon hole argument that
\begin{eqnarray*}
&&\forall j\in J \ \exists\tilde j(j)\in \tilde J\ \exists \text{ pairwise
distinct }\beta^{(1)}(j),\ldots,\beta^{(M)}(j):\\
&&\forall \eta=1,\ldots,M: \tilde j(j,\beta^{(\eta)}(j))=\tilde
j(j,\beta^{(1)}(j)).\label{pigeon1}
\end{eqnarray*}
Analogously, it follows from \eqref{bbb} and $\#J=\delta_0$ that
\begin{eqnarray*}
&&\forall \tilde j\in \tilde J \ \exists j(\tilde j)\in J\ \exists \text{
pairwise distinct }\beta^{(1)}(\tilde j),\ldots,\beta^{(M)}(\tilde j):\\
&&\forall \eta=1,\ldots,M: j(\tilde j,\beta^{(\eta)}(\tilde j))=j(\tilde
j,\beta^{(1)}(\tilde j)).\label{pigeon2}
\end{eqnarray*}
Because of the nondegeneracy condition \eqref{nondeg}, this implies that
\begin{equation*}
\forall j\in J\ \exists \tilde j\in \tilde J\ \forall \alpha=1,\ldots,M: g_\alpha^{(j)}(x)\equiv  \tilde g_\alpha^{(\tilde j)}(x).
\end{equation*}
By symmetry one also has that
\begin{equation*}
\forall \tilde j\in \tilde J\ \exists j\in J\ \forall \alpha=1,\ldots,M:
g_\alpha^{(j)}(x)\equiv \tilde g_\alpha^{(\tilde j)}(x).
\end{equation*}
The Identity Theorem now gives $C=\tilde C$.
\end{proof}
We can now go forward with the construction of our crucial embedding. Let $l:=4m-3$ (so that Proposition \ref{Clgen} can be applied to our
case). Let
\begin{equation*}
u:\P^m\xrightarrow{|lH|}\P(H^0(\P^m,lH))\cong\P^{{l+m\choose m}-1}
\end{equation*}
be the Veronese embedding of degree $l$ (as usual given by monomials $Z^\mu$
of multidegree $\mu$ with $|\mu|=l$). Let
\begin{equation*}
M:={l+m\choose m}-1.
\end{equation*}
Moreover, let
\begin{equation*}
\pi_\nu:\P^M\backslash\P^{M-3}_\nu\to\P^2_\nu\quad (\nu=1,\ldots,Q)
\end{equation*}
be the projections from Lemma \ref{injlin} with $\delta_0:=ld$, i.e.
\begin{equation*}
Q:=ldM{4ld+M\choose M}.
\end{equation*}\par
We take
\begin{eqnarray*}
\tilde u &:& \P^m\to\P^M\\
&&[Z]\mapsto[(\tilde f_\mu(Z))_\mu]
\end{eqnarray*}
to be a generic slight perturbation of $u$, chosen so that $\tilde f_\mu(Z)$
is a generic polynomial of degree $l$ close to the monomial $Z^\mu$ and such
that $\tilde u$ has the following properties. First, $\tilde u$ is still
injective. (This is clearly true for every small perturbation.) Second, for
every curve of degree $d$ in $\P^m$, $\overline{\pi_\nu(\tilde u(C)\cap\Dom(\pi_\nu))}$ is a
curve of degree $ld$ in $\P^2_\nu$ (when considered with the appropriate multiplicity). A choice of the $\tilde f_\mu$ with this
property is indeed possible. To see this, note that $\pi_\nu$ is given by
\begin{eqnarray*}
\pi_\nu&:&\P^M\backslash\P^{M-3}_\nu\to\P^2_\nu\\
&&[(Z_\mu)]\mapsto [\sum_{|\mu|=l}a^{(\nu,1)}_\mu
Z_\mu,\sum_{|\mu|=l}a^{(\nu,2)}_\mu Z_\mu,\sum_{|\mu|=l}a^{(\nu,3)}_\mu
Z_\mu]
\end{eqnarray*}
with $a^{(\nu,i)}_\mu\in \C$. Thus, we can apply Proposition \ref{Clgen} to
the triple
\begin{equation*}
f^{(\nu)}_i:=\sum_{|\mu|=l}a^{(\nu,i)}_\mu \tilde f_\mu \quad (i=1,2,3)
\end{equation*}
once the $\tilde f_\mu$ are chosen such that the triples
$f^{(\nu)}_1,f^{(\nu)}_2,f^{(\nu)}_3$ are generic in the sense of
Proposition \ref{Clgen}.\par
Summarizing, we have established the existence of rational maps
\begin{equation*}
\phi_\nu:=\pi_\nu\circ\tilde u:\P^m\rightharpoonup \P^2_\nu
\end{equation*}
such that the induced map
\begin{equation*}
\Phi_\nu:\Chow_{1,d}(\P^m)\to\Chow_{1,ld}(\P^2_\nu)\cong\P^{{ld+2\choose
2}-1}_\nu
\end{equation*}
given by
\begin{equation*}
[C]\mapsto [\overline{\phi_\nu(C\cap \tilde u^{-1}(\Dom(\pi_\nu)))}]
\end{equation*}
 (with appropriate multiplicity) is holomorphic. The map
\begin{eqnarray*}
\Phi&:&\Chow_{1,d}(\P^m)\to\prod_{\nu=1}^Q\Chow_{1,ld}(\P^2_\nu)\\
&&[C]\mapsto (\Phi_\nu([C]))_{\nu=1,\ldots,Q}
\end{eqnarray*}
is injective due to the injectivity of $\tilde u$ and Lemma \ref{injlin}, giving us the sought after holomorphic injection of $\Chow_{1,d}(\P^m)$ into a product of projective spaces. We shall refer to $\Phi$ as an embedding, although, strictly speaking, we have only proven it to be holomorphic and injective. Since we are only interested in counting degrees, these properties are all that we need for our arguments.\par
We conclude this subsection by establishing a bound for the total degree of the defining equations of $\Phi(\Chow_{1,d}'(\P^m))$. Here, $\Chow_{1,d}'(\P^m)$ is supposed to denote the union of those irreducible components of $\Chow_{1,d}(\P^m)$ whose general members are irreducible cycles. Clearly, we have that
\begin{equation*}
\psi_X(B)\subset\Chow_{1,d}'(\P^m)\subset \Chow_{1,d}(\P^m).
\end{equation*}
So far, this fact has been irrelevant and has thus been disregarded for generality's sake, but in the subsequent estimates, it makes things a little less involved.\par
From now on, we let $[A_\alpha^{(\nu)}]$ be homogeneous coordinates on $\Chow_{1,ld}(\P^2_\nu)\cong\P^{{ld+2\choose
2}-1}_\nu$, indexed by a multi-index $\alpha=(\alpha_0,\alpha_1,\alpha_2)$ with $|\alpha|=\alpha_0+\alpha_1+\alpha_2=ld$. We now consider the inclusion
\begin{gather*}
\Phi(\Chow_{1,d}'(\P^m))\subset \mathcal C := \left\{([A_\alpha^{(1)}],\ldots,[A_\alpha^{(Q)}])\in \prod_{\nu=1}^Q\Chow_{1,ld}(\P^2_\nu)|\right.\\
 \left.\bigcap_{\nu=1}^{Q}\{\sum_{|\alpha|=ld} A_\alpha^{(\nu)}(f^{(\nu)}_0)^{\alpha_0}(f^{(\nu)}_1)^{\alpha_1}(f^{(\nu)}_2)^{\alpha_2}=0\} \text{ has }\dim \geq 1 \text{ at all its points} \right\}.
\end{gather*}
Note that if we write $\mathcal C$ as the union $\bigcup_{j\in J}I_j$ of its irreducible components, then there exists a subset $J_1\subset J$ such that 
\begin{equation*}
\Phi(\Chow_{1,d}'(\P^m)) = \bigcup_{j\in J_1}I_j.
\end{equation*}
In order to prove this statement, one needs to ascertain that if $([A_\alpha^{(1)}],\ldots,[A_\alpha^{(Q)}])$ is a general point in $\Phi(\Chow_{1,d}'(\P^m))$, then $\bigcap_{\nu=1}^{Q}\{\sum_{|\alpha|=ld} A_\alpha^{(\nu)}(f^{(\nu)}_0)^{\alpha_0}(f^{(\nu)}_1)^{\alpha_1}(f^{(\nu)}_2)^{\alpha_2}=0\}$ is an irreducible curve and there exist small open sets $U_1$ in $\Chow_{1,d}'(\P^m)$ and $U_2\ni ([A_\alpha^{(1)}],\ldots,[A_\alpha^{(Q)}])$ in $\mathcal C$ such that $\Phi:U_1\to U_2$ is bijective. However, this is true because of the injectivity of $\Phi$ and the fact that for a general point in $\Phi(\Chow_{1,d}'(\P^m))$ we have that 
\begin{eqnarray*}
&&([A_\alpha^{(1)}],\ldots,[A_\alpha^{(Q)}])=\Phi([C])\\
&\Leftrightarrow & \bigcap_{\nu=1}^{Q}\{\sum_{|\alpha|=ld} A_\alpha^{(\nu)}(f^{(\nu)}_0)^{\alpha_0}(f^{(\nu)}_1)^{\alpha_1}(f^{(\nu)}_2)^{\alpha_2}=0\}=C,
\end{eqnarray*}
and, in addition, the ideal sheaf of $C$ is locally generated by $$ \sum_{|\alpha|=ld} A_\alpha^{(\nu)}(f^{(\nu)}_0)^{\alpha_0}(f^{(\nu)}_1)^{\alpha_1}(f^{(\nu)}_2)^{\alpha_2}\quad{\rm
for\ }1\leq\nu\leq Q.
$$
To achieve this additional condition, in our generic choice of $f^{(\nu)}_0,f^{(\nu)}_1,f^{(\nu)}_2$, we have to do the following (which we assume was already done at the time the choice was made). From each branch of $\Chow_{1,d}'(\P^m)$ choose one point which is represented by an irreducible cycle $C_\iota$ so that we get a finite collection $\left\{C_\iota\right\}_\iota$.  We choose $f^{(\nu)}_0,f^{(\nu)}_1,f^{(\nu)}_2$ so that for each $\iota$ the ideal sheaf of $C_\iota$ is locally generated by all $$ \sum_{|\alpha|=ld} A_\alpha^{(\nu)}(f^{(\nu)}_0)^{\alpha_0}(f^{(\nu)}_1)^{\alpha_1}(f^{(\nu)}_2)^{\alpha_2}
$$
which vanish identically on $C_\iota$.\par
We remark that the general points $([A_\alpha^{(1)}],\ldots,[A_\alpha^{(Q)}])$ of those irreducible components $I_j$ with $j\in J\backslash J_1$  have the property that $\bigcap_{\nu=1}^{Q}\{\sum_{|\alpha|=ld} A_\alpha^{(\nu)}(f^{(\nu)}_0)^{\alpha_0}(f^{(\nu)}_1)^{\alpha_1}(f^{(\nu)}_2)^{\alpha_2}=0\}$ is a reducible cycle. We also remark that the seeming contradiction in our argument stemming from the condition $\dim \geq 1$ (and not $\dim =1$) in the definition of $\mathcal C$ is in fact none. Namely, for certain special points in $\bigcup_{j\in J_1}I_j$ it is possible that there is a $[C]\in \Chow_{1,d}'(\P^m)$ such that if $([A_\alpha^{(1)}],\ldots,[A_\alpha^{(Q)}])=\Phi([C])$, then $\bigcap_{\nu=1}^{Q}\{\sum_{|\alpha|=ld} A_\alpha^{(\nu)}(f^{(\nu)}_0)^{\alpha_0}(f^{(\nu)}_1)^{\alpha_1}(f^{(\nu)}_2)^{\alpha_2}=0\}$ contains $C$ but has dimension greater than $1$.
\par
We now proceed to bounding the degree of defining equations for $\Phi(\Chow_{1,d}'(\P^m)) = \bigcup_{j\in J_1}I_j$ as follows. There is a Zariski open and dense set $U\subset  \Phi(\Chow_{1,d}'(\P^m))$ such that for $([A_\alpha^{(1)}],\ldots,[A_\alpha^{(Q)}])\in U$ we have the following. Let $h=\sum_{j=0}^{m}\xi_jZ_j$ be a linear form with indeterminate coefficients. 
\begin{eqnarray}
&&\dim (\bigcap_{\nu=1}^{Q}\{\sum_{|\alpha|=ld} A_\alpha^{(\nu)}(f^{(\nu)}_0)^{\alpha_0}(f^{(\nu)}_1)^{\alpha_1}(f^{(\nu)}_2)^{\alpha_2}=0\}) = 1\notag \\
&\Leftrightarrow&\dim (\bigcap_{\nu=1}^{Q}\{\sum_{|\alpha|=ld} A_\alpha^{(\nu)}(f^{(\nu)}_0)^{\alpha_0}(f^{(\nu)}_1)^{\alpha_1}(f^{(\nu)}_2)^{\alpha_2}=0\}) \geq 1\notag \\
&\Leftrightarrow&\dim \left(\left(\bigcap_{\nu=1}^{Q}\{\sum_{|\alpha|=ld} A_\alpha^{(\nu)}(f^{(\nu)}_0)^{\alpha_0}(f^{(\nu)}_1)^{\alpha_1}(f^{(\nu)}_2)^{\alpha_2} =0\}\right)\cap \{h=0\}\right) \geq 0,
\end{eqnarray}
where the set in the last line is to be understood as a subset of projective $m$-space over the function field $\C(\xi_0,\ldots,\xi_m)$. In his book \cite{Kollarbook}, Koll\'ar observed that the effective Nullstellensatz (e.g.\ in the form of \cite[7.4.4.3]{Kollarbook}) can be used to transform the statement of the last line into the effective statement that the ideal generated by
\begin{eqnarray*}
\sum_{|\alpha|=ld} A_\alpha^{(1)}(f^{(1)}_0)^{\alpha_0}(f^{(1)}_1)^{\alpha_1}(f^{(1)}_2)^{\alpha_2},\ldots, \sum_{|\alpha|=ld} A_\alpha^{(Q)}(f^{(Q)}_0)^{\alpha_0}(f^{(Q)}_1)^{\alpha_1}(f^{(Q)}_2)^{\alpha_2},\ h
\end{eqnarray*}
in $\C(\xi_0,\ldots,\xi_m)[Z_0,\ldots,Z_m]$ does not contain the power $(Z_0,\ldots,Z_m)^{(m+1)l^2d-m}$ of the irrelevant ideal. Koll\'ar suggests that this condition can be expressed as 
\begin{equation*}
\rank M < {(m+1)l^2d\choose m},
\end{equation*} 
where $M$ is a certain matrix whose nonzero entries are $\C$-linear combinations of the $A_\alpha^{(\nu)}$ (with fixed $\nu$) or one of the indeterminates $\xi_0,\ldots,\xi_m$.\par
We can now obtain a bound on the total degree of a set of defining equations for $\Phi(\Chow_{1,d}'(\P^m))$, as it follows from the minors-criterium for rank that $\Phi(\Chow_{1,d}'(\P^m))$ can be described by multi-homogeneous equations in the $A_\alpha^{(\nu)}$ whose total degree is less than ${(m+1)l^2d\choose m}$. The reader can find the missing details of this argument, which are left as an exercise in \cite{Kollarbook}, worked out in \cite[page 3]{Guerra}.\par
Finally, the subsequent Lemma \ref{segredegreebound} allows us to conclude that, after the Segre embedding
\begin{equation*}
\prod_{\nu=1}^Q\Chow_{1,ld}(\P^2_\nu)\hookrightarrow\P^{{ld+2\choose 2}^Q-1},
\end{equation*}
the variety $\Phi(\Chow_{1,d}'(\P^m))$ can be described by equations of degree also no more than 
\begin{equation*}
\max\left\{Q,{(m+1)l^2d\choose m}\right\}=Q.
\end{equation*}
\begin{lemma}\label{segredegreebound}
Let 
\begin{equation*}
\mathfrak s:\P^k\times\ldots\times \P^k\to\P^{(k+1)^r-1}
\end{equation*}
be the Segre embedding of the $r$-fold product of $\P^k$. Let a hypersurface $V$ of this product be defined by a multi-homogeneous equation $F$ of multidegree $(d_1,\ldots,d_r)$. Then $\mathfrak s(V)\subset \P^{(k+1)^r-1}$ is defined by equations of degree at most $\max \{r,d_1+\ldots+d_r\}$.
\end{lemma}
\begin{proof}
In homogeneous coordinates, $\mathfrak s$ is given by
\begin{equation*}
\mathfrak s ([X^{(1)}_0,\ldots,X^{(1)}_k],\ldots,[X^{(r)}_0,\ldots,X^{(r)}_k])=[(Z_{i_1,\ldots,i_r}:=X^{(1)}_{i_1}\cdot\ldots\cdot X^{(r)}_{i_r})_{i_1,\ldots,i_r=0,\ldots,k}].
\end{equation*}
Therefore, $\mathfrak s(V)$ is defined by equations of the form
\begin{equation*}
F_{i_1,\ldots,i_r}:=F(Z_{0,i_2,\ldots,i_r},\ldots,Z_{k,i_2,\ldots,i_r},\ldots,Z_{i_1,\ldots,i_{r-1},0},\ldots,Z_{i_1,\ldots,i_{r-1},k}),
\end{equation*} 
which are of degree $d_1+\ldots+d_r$, together with the defining equations of $\mathfrak s (\P^k\times\ldots\times \P^k)$, which are of degree $r$.
\end{proof}
\subsubsection{Obtaining degree bounds through multivalued sections}
In the sequel, we shall bound the degree of $\Phi_\nu(\psi_X(B))$, $\nu=1,\ldots,Q$. \par Let $N:=(ld)^2+1$ and, again, $l:=4m-3$ with $m:=1250(gq+s)$. According to \cite[Proposition 3]{P}, there exists a commutative base change diagram 
\begin{equation*} 
\begin{CD} 
X'@>\sigma_1>>X\times_BB'@>\sigma_2>>X\\
@Vf'VV@VVV @VVfV\\
B'@>{\id}>>B'@>\tau>>B\\
\end{CD}
\end{equation*}
with the following properties. $X'$ is the smooth minimal model of $X\times_BB'$. The family $f':X'\to B'$ has at least $N$ sections $s_i':B'\to X'$ $(i=1,\ldots,N)$ and $q'\leq 100N^2(gq+s)$ and $\deg \tau\leq 200N^2(gq+s)$. Moreover, the degrees with respect to $K_X$ of the multivalued sections $\sigma(s'_i(B'))$ ($\sigma:=\sigma_2\circ\sigma_1$) of the family $X\to B$ are bounded by 
\begin{equation*} 
\forall i=1,\ldots,N:\ \sigma(s'_i(B')).K_X\leq 100N(gq+s), 
\end{equation*} 
{\it i.e.,} 
\begin{equation}\label{bound} 
\forall i=1,\ldots,N:\ \deg \varphi_{|5K_X|}(\sigma(s'_i(B')))\leq  500N(gq+s). 
\end{equation} \par 
Fix $1\leq\nu\leq Q$ and choose a generic point $b_\nu'$ of $B'$ such that the $N$ points $(\phi_\nu\circ\varphi_{|5K_X|}\circ\sigma\circ
s_i')(b_\nu')$ ($1\leq i\leq N$) in $\P^2_\nu$ are distinct. \par 
By B\'ezout's Theorem an irreducible curve of degree $ld$ in the complex projective plane $\P^2_\nu$ is determined by any $N=(ld)^2+1$ pairwise distinct points on it (since any other irreducible curve of degree $ld$ sharing those points must coincide with the given curve). The Chow point $[A^{(\nu)}_\alpha(b_\nu')] \in \Chow_{1,ld}(\P^2_\nu)\cong\P^{{ld+2\choose 2}-1}_\nu$ of $\Phi_\nu([\varphi_{|5K_X|}(X_{\tau(b_\nu')})])$, as a set of ${ld+2\choose 2}$ unknowns, is determined by the following system of $N$ linear equations: 
\begin{eqnarray}\label{system} 
&&\sum_{|\alpha|=ld}A^{(\nu)}_\alpha(b_\nu')(f_0^{(\nu)})^{\alpha_0}(f_1^{(\nu)})^{\alpha_1} (f_2^{(\nu)})^{\alpha_2}( (\phi_\nu\circ\varphi_{|5K_X|}\circ\sigma\circ s_1')(b_\nu'))=0\nonumber\\
&&\hskip 4cm \vdots\\ 
&&\sum_{|\alpha|=ld}A^{(\nu)}_\alpha(b_\nu')(f_0^{(\nu)})^{\alpha_0}(f_1^{(\nu)})^{\alpha_1}(f_2^{(\nu)})^{\alpha_2}( (\phi_\nu\circ\varphi_{|5K_X|}\circ\sigma\circ s_N')(b_\nu'))=0.\nonumber
\end{eqnarray}
From \eqref{system} we can select a subsystem of ${ld+2\choose 2}-1$ linear equations
\begin{eqnarray}\label{subsystem} 
&&\sum_{|\alpha|=ld}A^{(\nu)}_\alpha(b_\nu')(f_0^{(\nu)})^{\alpha_0}(f_1^{(\nu)})^{\alpha_1} (f_2^{(\nu)})^{\alpha_2}( (\phi_\nu\circ\varphi_{|5K_X|}\circ\sigma\circ s_{i_1}')(b_\nu'))=0\nonumber\\
&&\hskip 4cm \vdots\\ 
&&\sum_{|\alpha|=ld}A^{(\nu)}_\alpha(b_\nu')(f_0^{(\nu)})^{\alpha_0}(f_1^{(\nu)})^{\alpha_1}(f_2^{(\nu)})^{\alpha_2}( (\phi_\nu\circ\varphi_{|5K_X|}\circ\sigma\circ s_{i_{{ld+2\choose 2}-1}}')(b_\nu'))=0.\nonumber
\end{eqnarray}
such that, for some $\alpha^{(0)}$, the $\left({ld+2\choose 2}-1\right)\times\left({ld+2\choose 2}-1\right)$ determinant of the matrix of coefficients obtained after removing the column corresponding to the variable $A^{(\nu)}_{\alpha^{(0)}}(b_\nu')$ is nonzero. \par 
Consider the following system of linear equations obtained from \eqref{subsystem} by replacing the point $b_\nu'$ of $B'$ by a variable point $b'$ of $B'$. 
\begin{eqnarray}\label{subsystem_at_variable_point}
&&\sum_{|\alpha|=ld}A^{(\nu)}_\alpha(b')(f_0^{(\nu)})^{\alpha_0}(f_1^{(\nu)})^{\alpha_1}(f_2^{(\nu)})^{\alpha_2}( (\phi_\nu\circ\varphi_{|5K_X|}\circ\sigma\circ s_{i_1}')(b'))=0\nonumber\\
&&\hskip 4cm \vdots\\ 
&&\sum_{|\alpha|=ld}A^{(\nu)}_\alpha(b')(f_0^{(\nu)})^{\alpha_0}(f_1^{(\nu)})^{\alpha_1}(f_2^{(\nu)})^{\alpha_2}( (\phi_\nu\circ\varphi_{|5K_X|}\circ\sigma\circ s_{i_{{ld+2\choose 2}-1}}')(b'))=0.\nonumber
\end{eqnarray}
We arrange the set of all $\alpha$ with $|\alpha|=ld$ so that the position of $\alpha$ is the integer $k_\alpha$ with $1\leq k_\alpha\leq{ld+2\choose 2}$. For every $\alpha$ with $|\alpha|=ld$ let $D_\alpha(b')$ be $(-1)^{k_\alpha-1}$ times the $\left({ld+2\choose 2}-1\right)\times\left({ld+2\choose 2}-1\right)$ determinant of the coefficients of the system \eqref{subsystem_at_variable_point} in the set of ${ld+2\choose 2}$ unknowns $(A^{(\nu)}_\beta(b'))_{|\beta|=ld}$ after the column corresponding to the variable $A^{(\nu)}_\alpha(b')$ is removed. \par
Choose two nonzero ${ld+2\choose 2}$-tuples $(c_\alpha)$ and $(d_\alpha)$ of complex numbers such that the intersection of the two hyperplanes in $\P^{{ld+2\choose 2}-1}_\nu$ defined by them is disjoint from the image of $B'$ under $\Phi_\nu\circ\psi_X\circ\tau$ and both $\sum_{|\alpha|=ld}c_\alpha D_\alpha(b_\nu')$ and $\sum_{|\alpha|=ld}d_\alpha D_\alpha(b_\nu')$ are nonzero. It follows from Cramer's rule applied to \eqref{subsystem_at_variable_point} that
\begin{equation*}
\frac{\sum_{|\alpha|=ld} c_\alpha A^{(\nu)}_\alpha(b')}{\sum_{|\alpha|=ld}d_\alpha A^{(\nu)}_\alpha(b')}=\frac{\sum_{|\alpha|=ld}c_\alpha D_\alpha(b')}{\sum_{|\alpha|=ld}d_\alpha D_\alpha(b')}
\end{equation*}
for all $b'\in B'$. \par 
Since $(f_0^{(\nu)})^{\alpha_0}(f_1^{(\nu)})^{\alpha_1} (f_2^{(\nu)})^{\alpha_2}((\phi_\nu\circ\varphi_{|5K_X|}\circ\sigma\circ s_i')(b'))$, as $b'$ varies as a point of $B'$, is a holomorphic section of $(s_i')^*(\sigma^*(\varphi_{|5K_X|} ^*({\mathcal O}_{\P^m}(l^2d)))) $ over $B'$, it follows that $D_i(b')$, as $b'$ varies as a point of $B'$, is a holomorphic section of 
\begin{eqnarray*} 
F&:=&\det \bigoplus_{j=1}^{{ld+2\choose 2}-1}(s_{i_j}')^*(\sigma^*(\varphi_{|5K_X|}^*({\mathcal O}_{\P^m}(l^2d))))\\
&\cong&\bigotimes_{j=1}^{{ld+2\choose 2}-1}(s_{i_j}')^* (\sigma^*(\varphi_{|5K_X|}^*({\mathcal O}_{\P^m}(l^2d))))
\end{eqnarray*}
over $B'$.
Now note that $\deg F$ and therefore $\deg(\Phi_\nu(\psi_X(B)))$ are no more than 
\begin{eqnarray*} 
&&\sum_{j=1}^{{ld+2\choose 2}-1}\int_{\sigma (s_{i_j}'(B'))} c_1 (\varphi_{|5K_X|}^*({\mathcal O}_{\P^m}(l^2d)) )\\
&\stackrel{\eqref{bound}}{\leq}&\left({l d+2\choose 2}-1\right)\cdot l^2d\cdot 500N(gq+s), 
\end{eqnarray*} 
which means we have achieved our goal of effectively bounding $\deg(\Phi_\nu(\psi_X(B)))$. To shorten notation, we denote the above expression by $D$, i.e. 
\begin{equation*} 
D:=\left({l\cdot 5(2g-2)+2\choose 2}-1\right) \cdot l^2\cdot 5(2g-2)\cdot 500((l\cdot 5(2g-2))^2+1)(gq+s). 
\end{equation*} 
\subsubsection{Bounding the degree of the image of the moduli map}
Now that we have established an effective bound for
the degree of $\Phi_\nu(\psi_X(B))$, we can easily bound the degree of $\Phi(\psi_X(B))$
in $\P^{{ld+2\choose 2}^Q-1}$. Namely, we have
\begin{eqnarray*}
&&\deg \Phi(\psi_X(B))\\
&=&\int_{\Phi(\psi_X(B))}\sum_{\nu=1}^Q{\rm pr}_\nu^*(\omega_{\text{Fubini-Study}})\\
&=&\sum_{\nu=1}^Q\deg \Phi_\nu(\psi_X(B))\\
&\leq& Q\cdot D.
\end{eqnarray*}
\subsection{Proof of the uniform effective bound}
We are now in a position to give a proof of Theorem \ref{effshafconj}.
\begin{proof}[Proof of Theorem \ref{effshafconj}]
First, take $q\geq 2$. We have proved that the degree of
\begin{equation*}
\Gamma_{\psi_X}\hookrightarrow
B\times\Phi(\Chow'_{1,d}(\P^m))\hookrightarrow\P^{5(q-1)A-1}
\end{equation*}
with
\begin{equation*}
A:={ld+2\choose 2}^Q
\end{equation*}
is bounded by
\begin{equation*}
\deg(3K_B)+\deg(\Phi(\psi_X(B)))\leq 6(q-1)+Q\cdot D.
\end{equation*}
Let $I$ be an irreducible component of ${\Chow}_{1,\tilde d}'(B\times\Phi(\Chow_{1,d}'(\P^m)))$ that contains the graph of one of our moduli maps. The arguments used in the proof of Parshin's original parametrization statement \cite[Theorem 1]{P} also apply to our situation. They tell us that there is a Zariski-open subset $U\subset I$ such that all Chow points $[\Gamma]\in U$ correspond to smooth families $V_{\Gamma}\to B\backslash S$. Moreover, we can infer from the rigidity statement proved by Arakelov (\cite[Proposition 2.1]{A}) that for any $[\Gamma_1],[\Gamma_2]\in U$ the corresponding families $V_{\Gamma_1}\to B\backslash S$ and $V_{\Gamma_2}\to B\backslash S$ are isomorphic. Therefore, the number of isomorphism classes of families in Theorem \ref{effshafconj} is no greater than the sum of the numbers of irreducible components of all the Chow
varieties
\begin{equation*}
{\Chow}_{1,\tilde d}'(B\times\Phi(\Chow_{1,d}'(\P^m))),
\end{equation*}
for $\tilde d=1,\ldots,6(q-1)+Q\cdot D$. However, from \cite{Kollarbook} and \cite{Guerra} the following Proposition is known.
\begin{proposition}
Let $W\subset \P^n$ be a projective variety defined by equations of degree no more than $\delta_1$. Let $\Chow'_{k,\delta_2}(W)$ denote the union of those irreducible components of $\Chow_{k,\delta_2}(W)$ whose general element is irreducible. Then the number of irreducible components of $\Chow'_{k,\delta_2}(W)$ is no more than
\begin{equation*}
{(n+1)\max\{\delta_1,\delta_2\}\choose n}^{(n+1)(\delta_2 {\delta_2+k-1\choose k }+{\delta_2+k-1 \choose k-1})}.
\end{equation*} 
\end{proposition}
Since it was shown at the end of Section \ref{twistedproj} that the degree of the defining functions of $B\times \Phi(\Chow_{1,d}'(\P^m))$ is no more than $Q$ ($< 6(q-1)+ Q\cdot D$), we can conclude that the number of isomorphism classes of families in Theorem \ref{effshafconj} is no greater than
\begin{equation*}
(6(q-1)+Q\cdot D) \cdot {5(q-1)A\cdot (6(q-1)+Q\cdot D))\choose 5(q-1)A-1}^{5(q-1)A((6(q-1)+Q\cdot D)^2+1)}.
\end{equation*}
This proves the case $q\geq 2$.\par
To handle the remaining two cases $q=0,1$, one is naturally led to making a
base change to a base curve $B'$ of genus $2$. For $q=0$, i.e. $B=\P^1$,
take any (hyperelliptic) curve $B'$ of genus two with its natural degree $2$
map (branched in $6$ points) to $\P^1$ as the base change map. The bad set
$S'$ will have cardinality at most $2s$, and it is known that the number of
families $\tilde X\to B$ such that $\tilde X\to B$ is not isomorphic to
$X\to B$, but $\tilde X\times_BB'$ is isomorphic to $X\times_B B'\to B'$, is
no more than the number of maps in the de Franchis-Severi
Theorem with a curve of genus $g$ as domain (independently of the (fixed) $B$ and $B'$). A proof of this statement is
contained, for example, in \cite[Proposition 1.1]{C}. Thus, our bound is
$\mathcal S(g)$ times the bound from the case $q'=2,s'=2s$.\par
For $q=1$ the usual branched covering trick allows a base change to a curve
of genus $2$ branched over $2$ points, and the bound we seek is again
$\mathcal S(g)$ times the bound for the case $q'=2,s'=2s$.
\end{proof}
\subsection{Some concluding remarks on the proof}\label{finalremark}
It seems worthwhile to remark that we resorted to the
generalization of Clemens' theorem only to make the maps $\Phi_\nu$
holomorphic. In this way, we did not have to handle blow-ups when
determining the degree of $\Phi_\nu(\psi_X(B))$. However, it should also be possible
to use rational projections and deal with this problem directly.\par
The reason why we confined ourselves to the case of the fibers being curves
is the following. First, if the fiber dimension is at least two, not much
seems to be known as to the right kind of rigidity conditions. Mere nonisotriviality is no longer sufficient for rigidity and finiteness, as the example of products of hyperbolic curves as fibers clearly shows. Secondly, in our proof of the
boundedness part, we made use of B\'ezout's Theorem saying that the image of
a fiber in one of the $\P^2_\nu$ is completely determined by an effectively
finite number of points on it. For a variety of dimension at
least two, this is the case only if those points are in general position,
i.e. the multivalued cross-sections must be chosen such that their
intersections with the fibers are in general position after the projections.
Trying to overcome these difficulties will be a subject of the author's
research in the future.
\section{A uniform effective solution to the Mordell Conjecture}\label{effmanin}
\subsection{Statement of the bound on rational points}
In \cite{P}, Parshin shows how, through an argument now known as ``Parshin's
construction'' or ``Parshin's trick'', Theorem \ref{shaf} implies the
following statement, originally known as the Mordell Conjecture for function
fields and proved by Manin in \cite{M} (see also \cite{Grauert}).
\begin{theorem}
Let $X(\mathbb C (B))$ be a nonisotrivial curve of genus at least $2$
defined over $\mathbb C (B)$. Then the number of $\mathbb C (B)$-rational
points on $X(\mathbb C (B))$ is finite.
\end{theorem}
Recall that {\it nonisotrivial} in this case means that the canonical
minimal family $X\to B$ pertaining to $X(\mathbb C (B))$ is
nonisotrivial.\par
In this section, we shall demonstrate how our Theorem \ref{effshafconj} implies the following uniform effective version of the
Mordell Conjecture for function fields via Parshin's method. The result of proving such a bound is not new (e.g., see
\cite{Miyaoka}), but to derive such a bound via Parshin's previously ineffective trick nevertheless seems to be
an interesting application of our Theorem \ref{effshafconj}.
\begin{theorem}\label{effmordellconj}
Let $X(\mathbb C (B))$ be a nonisotrivial curve of genus $g\geq 2$ defined
over $\mathbb C (B)$. Let the canonical minimal family pertaining to
$X(\mathbb C (B))$ have no more than $s$ singular fibers. Then the number of
$\mathbb C (B)$-rational points on $X(\mathbb C (B))$ is no more than
\begin{eqnarray*}
&&\mathcal S(2+2^{2g+1}(g-1))\\
&\cdot&P(2+2^{2g+1}(g-1),C(g,q,s),2^{2g}(2^{2g}-1)\cdot
2^{2(1+2^{2g}(g-1))}\cdot 2s)\\
&\cdot&(2+2^{2g+1}(g-1))\\
&\cdot& (2q+s)\cdot (C(g,q,s)+1) !,
\end{eqnarray*}
where $P(g',q',s')$ is the effective bound proved in Theorem
\ref{effshafconj}, $C(g,q,s)$ is
\begin{equation*}
1+ 2^{2g}(2^{2g}-1)\cdot 2^{2(1+2^{2g}(g-1))}\cdot
2(q-1)+(2^{2g}(2^{2g}-1)\cdot 2^{2(1+2^{2g}(g-1))}\cdot 2-1) s
\end{equation*}
and $\mathcal S(g)$ is the bound for the effective de Franchis-Severi
Theorem from page {\rm \pageref{HSbound}}.
\end{theorem}
\subsection{Parshin's trick}
There exist numerous expositions dealing with Parshin's trick, which works
in both the function field and the number field case. Thus, our presentation
of the construction will be very concise and essentially serves the purpose
of introducing the relevant notation and making the present work self-contained for the reader's convenience. We proceed along the lines of \cite[Section 4]{C}. See also \cite{Szpiro}.\par
The key idea at the heart of Parshin's construction is the observation that
rational points of $X(\mathbb C (B))$ are in bijective correspondence to
sections of $X\to B$, the canonical minimal family obtained from
$X(\mathbb C (B))$. It seems reasonable to expect that one might be
able to get a handle on the number of these sections. We will
indeed do so by associating a family $X'\to B'$ to every section.
The construction is such that the number of eligible families $X'\to B'$ can be
estimated by means of Theorem \ref{effshafconj} and such that $X'\to B'$
defines the section that gave rise to it up to an effectively finite number
of possibilities.
\par
To start the construction, let $P$ be a $\mathbb C (B)$-rational point of
$X(\mathbb C (B))$. Let it correspond to the section $\sigma$ of the
pertaining canonical minimal family $f:X\to B$. Let $\Sigma:=\sigma(B)$. Our
goal is to construct a finite covering $\rho :B'\to B$ (ramified only over
$S$) and a fibration $X'\to B'$ such that every fiber $X'_{b'}$ is a finite
covering of $X_{\rho(b')}$, ramified only over $\sigma (\rho(b'))$.\par
The section $\sigma$ gives rise to a map $u:X\to\pic X/B$ by setting
$u(x):=x-\sigma(f(x))$. The multiplication by $2$ map on $\pic X/B$ yields a
covering of $X$ over $B$ that is \'etale of degree $2^{2g}$ outside of the
singular fibers of $X$. Let $Y$ be a connected component of this covering.
By the Riemann-Hurwitz Theorem, $Y\to B$ is a family of curves of genus at
most $1+2^{2g}(g-1)$. Denote the preimage of $\Sigma$ in $Y$ by $D$.\par
After a base change $B_1\to B$ of degree at most $2^{2g}(2^{2g}-1)$,
ramified only over $S$, there are two disjoint sections
$\sigma_1,\sigma_2:B_1\to Y_1$ (that both map to $D$) on the minimal
resolution $Y_1$ of $B_1\times_B Y$. Let $\Sigma_i=\sigma_i(B_1)$ for
$i=1,2$.\par
Next, we map $B_1$ to $\pic Y_1/B_1$ by setting
$u_1(b_1):=\sigma_1(b_1)-\sigma_2(b_1)$. We take $B_2$ to be a connected
component of the fibered product $B_1\times_{\pic Y_1/B_1}\pic Y_1/B_1$
under $u_1$ and under the multiplication by $2$ map $\pic Y_1/B_1\to\pic
Y_1/B_1$. On $Y_2:=Y_1\times_{B_1}B_2$, the line bundle $\mathcal
O_{Y_2}(\Gamma_1+\Gamma_2)$, defined as the pull-back line bundle of
$\mathcal O_{Y_1}(\Sigma_1+\Sigma_2)$, admits a square root. This ensures
that, after a degree 2 base change $B_3\to B_2$, there exists a double
covering $Y_3\to Y_2$ having branch locus $\Gamma_1+\Gamma_2$. After taking
$X'$ to be the minimal resolution of $Y_3$ and $B'$ to be $B_2$, we have met
our objective, as we shall see in the subsequent section.\par
\subsection{Proof of the bound on rational points}
Clearly, what we need to do is to bound the number of those families $X'\to
B'$ that can possibly occur from the Parshin construction and to bound the
number of those sections that may give rise to the same given family $X'\to
B'$.\par
The first thing to notice is that according to the Riemann-Hurwitz Theorem
\begin{equation*}
g'=\genus (X')\leq 2+2^{2g+1}(g-1),
\end{equation*}
because, for $b'\in B'\backslash\rho^{-1}(S)$, $X'_{b'}$ is a covering of
$X_{\rho(b')}$ of degree at most $2^{2g+1}$, with only two simple
ramification points over $\Sigma$.\par
Furthermore, the degree of $\rho:B'\to B$ is no more than
\begin{equation*}
2^{2g}(2^{2g}-1)\cdot 2^{2(1+2^{2g}(g-1))}\cdot 2,
\end{equation*}
because
\begin{itemize}
\item $\deg(B_1\to B)\leq 2^{2g}(2^{2g}-1)$
\item $\deg(B_2\to B_1)\leq 2^{2(1+2^{2g}(g-1))}$
\item $\deg(B_3\to B_2) = 2$.
\end{itemize}
By Riemann-Hurwitz, we find that $q'$ is no more than
\begin{equation*}
1+ 2^{2g}(2^{2g}-1)\cdot 2^{2(1+2^{2g}(g-1))}\cdot
2(q-1)+(2^{2g}(2^{2g}-1)\cdot 2^{2(1+2^{2g}(g-1))}\cdot 2-1) s.
\end{equation*}
Let $C(g,q,s)$ denote the above number.\par Next, we determine the number of
possibilities for $B'$. For this, we need the following Lemma.
\begin{lemma}\label{coverbound}
For any integer $\theta$, the number of isomorphism classes of smooth
Riemann surfaces $B'$ that allow a surjective holomorphic map of degree
$\theta$ to $B$ that is branched only over $S$ is at most
\begin{equation*}
(2q+s)\cdot\theta !.
\end{equation*}
\end{lemma}
\begin{proof}
It is well known that the fundamental group of $B\backslash S$ is generated
by $2q+s$ elements (with one relation among them). Thus, the statement we
seek to prove is immediate from the fact that such a $B'$ is determined up
to isomorphism by a homomorphism from the fundamental group of $B\backslash
S$ to the symmetric group on $\theta$ letters.
\end{proof}
According to Lemma \ref{coverbound}, the number of possible $B'$ can be
bounded by
\begin{eqnarray*}
&&\sum_{\theta=1}^{C(g,q,s)} (2q+s)\theta !\\
&\leq& C(g,q,s)\cdot(2q+s)\cdot C(g,q,s) !\\
&\leq&(2q+s)\cdot (C(g,q,s)+1) !.
\end{eqnarray*}
Now, if $P(g,q,s)$ is the effective bound for the Shafarevich Conjecture
proved in Theorem \ref{effshafconj}, then the number of families $X'\to B'$
that may arise from the Parshin construction is no more than
\begin{eqnarray*}
&&P(2+2^{2g+1}(g-1),C(g,q,s),2^{2g}(2^{2g}-1)\cdot 2^{2(1+2^{2g}(g-1))}\cdot
2s)\\
&\cdot&(2+2^{2g+1}(g-1))\\
&\cdot& (2q+s)\cdot (C(g,q,s)+1) !.
\end{eqnarray*}\par
To finish the proof of Theorem \ref{effmordellconj}, we need to estimate the
number of rational sections giving rise to the same family $X'\to B'$.\par
First, note that Parshin's construction is set up such that $\Sigma$ is
precisely the branch locus of the map $X'\to X$. Thus, the number of
rational sections yielding the same family $X'\to B'$ can be no more than
the number of $B$-maps $X'\to X$. Such maps, however, can be regarded as
maps $X'(\mathbb C (B'))\to X(\mathbb C (B'))$ of curves over the function
field $\mathbb C (B')$, and it is known (see \cite[Proposition 1.1]{C}) that
the bound $\mathcal S(g)$ from the effective de Franchis-Severi Theorem over
$\mathbb C$ is valid also in the function field case. Therefore, no more
than $\mathcal S(2+2^{2g+1}(g-1))$ distinct rational points give rise to the
same family $X'\to B'$.\par
Summing up, the number of $\mathbb C (B)$-rational points on $X(\mathbb C
(B))$ is no more than
\begin{eqnarray*}
&&\mathcal S(2+2^{2g+1}(g-1))\\
&\cdot&P(2+2^{2g+1}(g-1),C(g,q,s),2^{2g}(2^{2g}-1)\cdot
2^{2(1+2^{2g}(g-1))}\cdot 2s)\\
&\cdot&(2+2^{2g+1}(g-1))\\
&\cdot& (2q+s)\cdot (C(g,q,s)+1) !,
\end{eqnarray*}
as it was asserted in Theorem \ref{effmordellconj}.\hfill$\qed$\par
\begin{acknowledgement}
The results obtained in this article first appeared in the author's Ph.D. Dissertation, which was presented in 2002 to the Department of Mathematics of the Ruhr-Universit\"at Bochum. It is a great pleasure to thank my advisors Prof.\ A.\ Huckleberry and Prof.\ Y.-T.\ Siu for their support. I am particularly indebted to Prof. Siu for numerous discussions on the subject of (effective) algebraic geometry and complex analysis, which took place during extended visits to the Mathematics Department of Harvard University and the Institute of Mathematical Research at Hong Kong University (hosted by Prof. N. Mok). Finally, it is with sincere gratitude that I acknowledge support through the Studienstiftung des Deutschen Volkes and the Schwerpunktprogramm ``Globale Methoden in der komplexen Geometrie'' of the Deutsche Forschungsgemeinschaft.
\end{acknowledgement}

\end{document}